\documentclass[11pt]{amsart}
\usepackage{amssymb, latexsym}
\theoremstyle{plain}
\newtheorem{theorem}{Theorem}

\newtheorem*{C}{Theorem C}

\newtheorem*{2'}{Theorem 2'}
\newtheorem*{3'}{Theorem 3'}

\theoremstyle{remark}

\newtheorem*{Remark 1}{Remark 1}
\newtheorem*{Remark 2}{Remark 2}
\newtheorem*{Remark 3}{Remark 3}
\newtheorem*{Remark 4}{Remark 4}

\DeclareRobustCommand{\stirling}{\genfrac\{\}{0pt}{}}

\numberwithin{equation}{section}

\begin{document}

\title [Permutation Fixing  Set, Touchard Polynomial, Cover of Multiset]
{Some Connections Between Cycles and  Permutations that Fix a Set  and
 Touchard Polynomials and Covers of Multisets }

\author{Ross G. Pinsky}

%\noindent  pinsky@math.technion.ac.il\ \ \ \ tel: 972-4-829-4083\ \ \  fax: 972-4-829-3388

\address{Department of Mathematics\\
Technion---Israel Institute of Technology\\
Haifa, 32000\\ Israel}
\email{ pinsky@math.technion.ac.il}

\urladdr{http://www.math.technion.ac.il/~pinsky/}

\subjclass[2000]{60C05, 05A05, 05A15} \keywords{ permutations that
fix a set, covers of multisets,  Touchard Polynomials, Dob\'inski's formula, Bell numbers, cycles in random permutations,
Ewens sampling distribution }
\date{}

\begin{abstract}

We present a new proof of a fundamental result concerning cycles of random permutations which
gives some intuition for the connection between Touchard polynomials and the Poisson distribution.
We also introduce a rather novel  permutation statistic
and study its distribution. This quantity, indexed by $m$, is the number of sets of size $m$ fixed by the permutation.
This leads  to a new
and simpler derivation of the exponential generating function for the number of covers of certain multisets.

\end{abstract}

\maketitle
\section{Introduction and Statement of Results}
In this paper, we present a new and simpler proof of a fundamental result concerning cycles of random permutations which
gives some intuition for the connection between Touchard polynomials and the Poisson distribution.
We also introduce a rather novel permutation statistic and study its distribution.
This quantity, indexed by $m$, is the number of sets of size $m$ fixed by the permutation.
This leads  to a new
and simpler derivation of the exponential generating function for the number of covers of certain multisets.

We begin by recalling some basic facts concerning Bell numbers and Touchard polynomials, and their
connection to Poisson distributions. The facts noted below without proof can be found in many books
on combinatorics; for example, in \cite{S}, \cite{W}.
The Bell number $B_n$ denotes the number of partitions of a set of $n$ distinct elements. Elementary combinatorial reasoning yields the recursive formula
\begin{equation}\label{Bell}
B_{n+1}=\sum_{k=0}^n\binom n k B_k,\ n\ge0,
\end{equation}
where $B_0=1$.
Let
\begin{equation}\label{Eseries}
E_B(x)=\sum_{n=0}^\infty \frac{B_n}{n!}x^n,
\end{equation}
denote the exponential moment generating function of $\{B_n\}_{n=0}^\infty$. Using \eqref{Bell} it is easy to show that
$E_B'(x)=e^xE_B(x)$, from which it follows that
\begin{equation}\label{E_B}
E_B(x)=e^{e^x-1}.
\end{equation}

A random variable $X$ has the Poisson distribution Pois$(\lambda)$, $\lambda>0$, if
$P(X=k)=e^{-\lambda}\frac{\lambda^k}{k!},\ k=0,1,\dots$.
Let  $M_\lambda(t)=Ee^{tX}$ denote the moment generating function of $X$, and let $\mu_{n;\lambda}=EX^n$ denote the $n$th moment of $X$.
Since $M_\lambda^{(n)}(0)=\mu_{n;\lambda}$, we have
\begin{equation}\label{Mseries}
M_\lambda(t)=\sum_{n=0}^\infty \frac{\mu_{n;\lambda}}{n!}t^n.
\end{equation}
However a direct calculation gives
\begin{equation}\label{MGF}
M(t)=\sum_{k=0}^\infty e^{tk}e^{-\lambda}\frac{\lambda^k}{k!}=e^{\lambda(e^t-1)}.
\end{equation}
From \eqref{Eseries}-\eqref{MGF}, it follows that the $n$th moment $\mu_{n;1}$ of a Pois$(1)$-distributed random variable satisfies
\begin{equation}\label{Bellpois}
\mu_{n;1}=B_n.
\end{equation}
Since
$$
\mu_{n;1}=EX^n=\sum_{k=0}^\infty (e^{-1}\frac1{k!})k^n,
$$
we conclude that
\begin{equation}\label{Dobinski}
B_n=\frac1e\sum_{k=0}^\infty \frac{k^n}{k!},
\end{equation}
which is known as \it Dob\'inski's formula\rm.

The Stirling number of the second kind $\stirling nk$ denotes
the number of partitions of a set of $n$ distinct elements into $k$ nonempty sets.
Elementary combinatorial reasoning yields the recursive formula
$$
\stirling{n+1}k=k\stirling nk+\stirling n{k-1}.
$$
From this, it is  not hard to derive the formula
\begin{equation}\label{stirpoly}
x^n=\sum_{j=0}^n\stirling nj(x)_j,
\end{equation}
where $(x)_j=x(x-1)\cdots (x-j+1)$ is the  falling factorial, and one defines $(x)_0=\stirling00=1$.
Now using \eqref{stirpoly} and the fact that $(k)_j=0$ for $j>k$, we can write the $n$th moment $\mu_{n;\lambda}$ of a Pois$(\lambda)$-distributed random variable as
\begin{equation}\label{Touchardmoment}
\begin{aligned}
&\mu_{n;\lambda}=\sum_{k=0}^\infty (e^{-\lambda}\frac{\lambda^k}{k!})k^n=\sum_{k=0}^\infty (e^{-\lambda}\frac{\lambda^k}{k!})\Big(\sum_{j=0}^n\stirling nj(k)_j\Big)=\\
&\sum_{k=0}^\infty (e^{-\lambda}\frac{\lambda^k}{k!})\Big(\sum_{j=0}^{k\wedge n}\stirling nj(k)_j\Big)=e^{-\lambda}\sum_{j=0}^n\stirling nj\sum_{k=j}^\infty\lambda^k\frac{(k)_j}{k!}=\\
&\sum_{j=0}^n\stirling nj\lambda^j.
\end{aligned}
\end{equation}

The Touchard polynomials $T_n(x)$, $n\ge0$, are defined by
$$
T_n(x)=\sum_{j=0}^n\stirling njx^j.
$$
Thus, \eqref{Touchardmoment} gives the formula
\begin{equation}\label{Touchardpois}
\mu_{n;\lambda}=T_n(\lambda).
\end{equation}
Since
$$
\mu_{n;\lambda}=\sum_{k=0}^\infty (e^{-\lambda}\frac{\lambda^k}{k!})k^n,
$$
we conclude from \eqref{Touchardpois} that
\begin{equation}\label{Dobinski-gen}
T_n(x)=e^{-x}\sum_{k=0}^\infty\frac{k^n}{k!}x^k.
\end{equation}
Since $T_n(1)=\sum_{j=0}^n\stirling nj=B_n$, the Dob\'inski formula \eqref{Dobinski} is contained in \eqref{Dobinski-gen}.

Let $S_n$ denote the set of permutations of $[n]=:\{1,\ldots, n\}$. For $\sigma\in S_n$, let $C^{(n)}_m(\sigma)$ denote the number of cycles of length $m$ in $\sigma$.
Let $P_n$ denote the uniform probability measure on $S_n$. We can now think of $\sigma\in S_n$ as random, and of $C^{(n)}_m$ as a random variable.
 Using generating function techniques and/or inclusion-exclusion formulas,
one can show that under $P_n$, the distribution of the random variable $C^{(n)}_m$ converges weakly to $Z_\frac1m$, where $Z_{\frac1m}$  has the Pois$(\frac1m)$-distribution;  equivalently;
\begin{equation}\label{wkconvm}
\lim_{n\to\infty}P_n(C^{(n)}_m=j)=e^{-\frac1m}\frac1{m^jj!},\ j=0,1\ldots.
\end{equation}

More generally, we consider the Ewens sampling distributions, $P_{n;\theta}$, $\theta>0$, on $S_n$
as follows. Let $N^{(n)}(\sigma)$ denote the number
of cycles in the permutation $\sigma\in S_n$,
and let $s(n,k)=|\{\sigma\in S^n:N^{(n)}(\sigma)=k\}|$ denote the number of permutations in $S_n$ with $k$ cycles.
It is known that the polynomial $\sum_{k=1}^n s(n,k)\theta^k$ is equal to the rising factorial
$\theta^{(n)}$, defined by $\theta^{(n)}=\theta(\theta+1)\cdots(\theta+n-1)$.
For $\theta>0$, define the probability measure
$P_{n;\theta}$ on $S_n$ by
$$
P_{n;\theta}(\{\sigma\})=\frac{\theta^{N^{(n)}(\sigma)}}{\theta^{(n)}}.
$$
Of course, $P_{n;1}$ reduces to the uniform measure $P_n$.
The following theorem can be proven; see for example, \cite{A}, \cite{P}.

\begin{C}\label{C}
Under $P_{n,\theta}$, the random vector $(C^{(n)}_1,C^{(n)}_2,\cdots ,C^{(n)}_m,\ldots)$ converges weakly to $(Z_{\theta},Z_\frac\theta2,\cdots, Z_{\frac\theta m},\cdots)$, where the random variables $\{Z_\frac\theta m\}_{m=1}^\infty$ are independent,
and $Z_\frac\theta m$ has the Pois$(\frac\theta m)$-distribution:
\begin{equation}\label{wkconvall}
(C^{(n)}_1,C^{(n)}_2,\cdots C^{(n)}_m,\ldots)\stackrel{w}{\Rightarrow}(Z_\theta,Z_\frac\theta 2,\cdots, Z_{\frac{\theta}m}\cdots);
\end{equation}
equivalently,
\begin{equation}\label{wkconvallexplicit}
\begin{aligned}
&\lim_{n\to\infty}P_n(C^{(n)}_1=j_1,C^{(n)}_2=j_2,\ldots, C^{(n)}_m=j_m)=\prod_{k=1}^me^{-\frac\theta k}\frac{(\frac{\theta}k)^{j_k}}{j_k!},\\
&\text{for all}\ m\ge1 \ \text{and}\ j_1,\ldots j_m\in \mathbb{Z}_+.
\end{aligned}
\end{equation}
\end{C}

We will use the method of moments to give a new and simpler proof of Theorem C, which will give
intuition for \eqref{Touchardpois}, or equivalently,
 for \eqref{Dobinski-gen}; that is for the connection between the moments of Poisson random variables and Touchard
 polynomials.

We now consider a permutation statistic that hasn't been studied much.
(Indeed, it was only after completing the first version of this paper that we were directed to any papers on this subject.)
For $\sigma\in S_n$ and  $A\subset[n]$,  define $\sigma(A)=\{\sigma_j:j\in A\}$.
If $\sigma(A)=A$, we will say that $\sigma$ fixes  $A$.
Let $\mathcal{E}^{(n)}_m(\sigma)$ denote the number of sets of cardinality $m$ that are fixed by $\sigma$. (Note that $\mathcal{E}^{(n)}_1(\sigma)=C^{(n)}_1(\sigma)$, the number of fixed points
of $\sigma$.) A little  thought reveals that
\begin{equation}\label{equ-en}
\mathcal{E}^{(n)}_m(\omega)=\sum_{\stackrel{(l_1,\ldots, l_m):\sum_{j=1}^mjl_j=m}{l_j\le C^{(n)}_j, j\in[m]}}\ \
\prod_{j=1}^m \binom{C^{(n)}_j(\omega)}{l_j}.
\end{equation}
For example, if $\sigma\in S_9$ is written in cycle notation as $\sigma=(379)(24)(16)(5)(8)$, then $\mathbb{E}^{(9)}_4(\omega)=5$, with the sets $A\subset[9]$ for which $|A|=4$ and
$\sigma(A)=A$ being  $\{3,5,7,9\}, \{3,7,8,9\}, \{1,2,4,6\}, \{2,4,5,8\}, \{1,5,6,8\}$.

We consider the uniform measure $P_n=P_{n;1}$ on $S_n$. From Theorem C and \eqref{equ-en} it follows that
the random variable $\mathcal{E}^{(n)}_m$ under $P_n$ converges weakly as $n\to\infty$ to the random variable
\begin{equation}\label{Einfty}
\mathcal{E}_m=:\sum_{\stackrel{(l_1,\ldots, l_m):\sum_{j=1}^mjl_j=m}{l_j\le Z_\frac1j, j\in[m]}}\ \ \prod_{j=1}^m \binom{Z_\frac1j}{l_j},
\end{equation}
where $\{Z_\frac1j\}_{j=1}^m$ are independent and $Z_\frac1j$ has the Pois$(\frac1j)$-distribution.
\medskip

\bf \noindent Remark.\rm\
Note that $\mathcal{E}_1=Z_1$, $\mathcal{E}_2=\binom{Z_1}2+Z_\frac12$,
$\mathcal{E}_3=Z_\frac13+Z_\frac12Z_1+\binom{Z_1}3$.
\medskip

For $k,m\in\mathbb{N}$, consider the multiset consisting of $m$ copies of the set $[k]$. A
 collection $\{\Gamma_l\}_{l=1}^r$ such that each $\Gamma_l$ is a nonempty subset of $[k]$, and such that each
 $j\in[k]$ appears in exactly $m$ from among the $r$ sets $\{\Gamma_l\}_{l=1}^r$,
is called an \it $m$-cover of $[k]$ of order $r$.\rm\
Denote the total number of $m$ covers of $[k]$, regardless of order, by
$v_{k;m}$. Note that when $m=1$, we have $v_{k;1}=B_k$, the $k$th Bell number, denoting the number of partitions of a set of $k$ elements.
Also, it's very easy to see that $v_{1;m}=1$ and $v_{2,m}=m+1$.

By calculating directly the moments of $\mathcal{E}^{(n)}_m$, we will prove the following theorem.
\begin{theorem}\label{1}
For $m,k\in\mathbb{N}$,
$$
E\mathcal{E}_m^k=v_{k;m}.
$$
In particular, $E\mathcal{E}_m=1$ and $E\mathcal{E}_m^2=m+1$; thus, Var$(\mathcal{E}_m)=m$.
\end{theorem}
%\bf\noindent Remark.\rm\ Consider the case $m=1$.
%As noted above,  $v_{k;1}=B_k$, the $k$th Bell number and   $e^{(n)}_1=C^{(n)}_1$.
%So by \eqref{wkconvm}, the random variable $X_1$ is distributed as Pois$(1)$, and Theorem \ref{1} reduces to
%\eqref{Bellpois}.

\bf \noindent Remark.\rm\ It is natural to suspect that $\mathcal{E}_m$ converges weakly to 0 as $m\to\infty$;
that is, $\lim_{n\to\infty}P(\mathcal{E}_m\ge1)=0$. This is in fact a hard problem.
In \cite{LP} it was shown that $P_n(\mathcal{E}_m^{(n)})\le Am^{-\frac1{100}}$, for $1\le m\le \frac n2$ and $n\ge2$.
Thus indeed, $\mathcal{E}_m$ converges weakly to 0 as $m\to\infty$. A lower bound on
$P(\mathcal{E}_m\ge1)$ of the form $A\frac{\log m}m$  was obtained in \cite{DFG}. These results were dramatically
improved in \cite{PPR} where it was shown that $P(\mathcal{E}_m\ge1)=m^{-\delta+o(1)}$ as $m\to\infty$,
where $\delta=1-\frac{1+\log\log2}{\log2}\approx0.08607$.
And very recently, in \cite{EFG},  this latter bound has been refined  to $A_1m^{-\delta}(1+\log m)^{-\frac32}\le P(\mathcal{E}_m\ge1)\le
A_2m^{-\delta}(1+\log m)^{-\frac32}$.
\medskip

Let
$$
V_m(x)=\sum_{k=1}^\infty \frac{v_{k;m}}{k!}x^k
$$
 denote the exponential generating function of the sequence $\{v_{k;m}\}_{k=1}^\infty$.
Of course, by Theorem \ref{1} $V_m$ is also the moment generating function of the random variable
$\mathcal{E}_m$: $V_m(x)=Ee^{x\mathcal{E}_m}$.
 Using \eqref{Einfty} and Theorem \ref{1}, we will give an almost immediate proof of the following representation theorem for $V_m(x)$. We use the notation $[z^m]P(z)=a_m$, where $P(z)=\sum_{m=0}^\infty a_mz^m$.
 \begin{theorem}\label{2}
 \begin{equation}
 V_m(x)=Ee^{x\mathcal{E}_m}=    e^{-\sum_{j=1}^m\frac1j}\sum_{u_1,\ldots, u_m\ge0}\big(\prod_{j=1}^m\frac{j^{-u_j}}{u_j!}\big)~e^{x\gamma_m(u)},
 \end{equation}
where
$$
\gamma_m(u)=[z^m]\prod_{j=1}^m(1+z^j)^{u_j}.
$$
 \end{theorem}
\noindent\bf Remark.\rm\ When $m=2,3$, the above formula reduces to
$$
\begin{aligned}
&V_2(x)=e^{-\frac32}e^{\frac12e^x}\sum_{r=0}^\infty \frac{e^{\binom r2x}}{r!};\\
&V_3(x)=e^{-\frac{11}6}e^{\frac{e^x}3}\sum_{r=0}^\infty\frac{e^{\binom r3x+\frac12e^{rx}}}{r!}.
\end{aligned}
$$
The formula for $m=2$ was proved by Comtets \cite{C} and the formula for $m=3$ was proved by Bender \cite{B}.
The case of general $m$ was proved by Devitt and Jackson \cite{DJ}.
They also prove  that there exists a number $c$ such that the extraction of the coefficient $v_{k;m}$ from the exponential generating function $V_m(x)$ can be done
in no more than $ck^m\log k$ arithmetic operations.

In section \ref{ProofC} we will give our new proof of Theorem C via the method of moments. In section \ref{Proofmulti}
we prove Theorems \ref{1} and \ref{2}.

\section{ A proof of Theorem C via the method of moments}\label{ProofC}
If a sequence of nonnegative random variables $\{X_n\}_{n=1}^\infty$ satisfies
$\sup_{n\ge1}EX_n<\infty$, then the sequence is tight, that is, pre-compact with respect to weak convergence.
Let $X$ be distributed as one of the accumulation points.
If for some $k\in \mathbb{N}$, $\lim_{n\to\infty}EX_n^k$ exists and equals $\mu_k$, and $\sup_{n\ge1}EX_n^{k+1}<\infty$,
then the $\{X_n^k\}_{n=1}^\infty$ are uniformly integrable, and thus $EX^k=\mu_k$. Thus,
if
\begin{equation}\label{moments}
\mu_k=:\lim_{n\to\infty}EX_n^k\ \text{exists for all}\ k\in\mathbb{N},
\end{equation}
 then $EX^k=\mu_k$, for all
$k$. The Stieltjes moment theorem states that if
\begin{equation}\label{Stie}
\sup_{k\ge1}\frac{\mu_k^\frac1k}k<\infty,
\end{equation}
 then
the sequence $\{\mu_k\}_{k=1}^\infty$ uniquely characterizes the distribution \cite{D}.
We conclude then that if a sequence of nonnegative random variables $\{X_n\}_{n=1}^\infty$ satisfies
\eqref{moments} and \eqref{Stie}, then the sequence is weakly convergent to a random variable $X$ satisfying
$EX^k=\mu_k$.

An extremely crude argument shows that the Bell numbers satisfy $B_k\le k^k$; thus
\begin{equation}\label{Bell-Stie}
\sup_{k\ge1}\frac{B_k^\frac1k}k<\infty.
\end{equation}
By \eqref{Touchardpois}, the $k$th moment $\mu_{k;\frac\theta m}$ of the Pois$(\frac\theta m)$-distributed random variable
$Z_\frac\theta m$ is equal to $T_k(\frac\theta m)$. Now $T_k(\frac\theta m)$ is  bounded from above by
 $T_k(\theta)$, for all $m\ge1$, and $T_k(\theta)\le B_k\max(1,\theta^k)$.
Thus, in light of \eqref{Bell-Stie} and the previous paragraph, if we prove
that
\begin{equation}\label{momentmethodCm}
\lim_{n\to\infty}E_{n;\theta}(C^{(n)}_m)^k=T_k(\frac\theta m),\ k,m\in\mathbb{N},
\end{equation}
where $E_{n;\theta}$ denotes the expectation with respect to $P_{n;\theta}$, then we will have proved that $C^{(n)}_m$ under $P_{n;\theta}$ converges weakly to $Z_\frac\theta m$, for all $m\in\mathbb{N}$.
And if we then prove  that
\begin{equation}\label{momentmethodCms}
\lim_{n\to\infty}E_{n;\theta}\prod_{j=1}^m(C^{(n)}_j)^{k_j}=\prod_{j=1}^mT_{k_j}(\frac\theta j),\ m\ge2, k_j\in\mathbb{N}, j=1,\ldots, m,
\end{equation}
then we will have completed the proof of Theorem C.

We first prove \eqref{momentmethodCm}.
In fact, we will first prove \eqref{momentmethodCm} in the case of the uniform measure, $P_n=P_{n;1}$.
Once we have this, the case of general $\theta$ will follow after a short explanation.
Assume that $n\ge mk$.
For $D\subset[n]$ with $|D|=m$,
let $1_D(\sigma)$ be equal to 1 or  0 according to whether or not $\sigma\in S_n$ possesses an $m$-cycle consisting of
the elements of $D$. Then we have
\begin{equation}\label{Cnrep}
C^{(n)}_m(\sigma)=\sum_{\stackrel{D\subset[n]}{|D|=m}}1_D(\sigma),
\end{equation}
and
\begin{equation}\label{Cnmkthmoment}
E_n(C^{(n)}_m)^k=\sum_{\stackrel{D_j\subset[n],|D_j|=m}{j\in[k]}}E_n\prod_{j=1}^k1_{D_j}.
\end{equation}

Now $E_n\prod_{j=1}^k1_{D_j}\neq0$ if and only if for some $l\in[k]$, there exist
disjoint sets $\{A_i\}_{i=1}^l$ such that $\{D_j\}_{j=1}^k=\{A_i\}_{i=1}^l$.
If this is the case, then
\begin{equation}\label{ldistinct}
E_n\prod_{j=1}^k1_{D_j}=\frac{(n-lm)!((m-1)!)^l}{n!}.
\end{equation}
(Here we have used the assumption that $n\ge mk$, since otherwise $n-ml$ will be negative for certain $l\in[k]$.)
The number of ways to construct $l$ disjoint, ordered sets  $\{A_i\}_{i=1}^l$, each of which consists of $m$ elements from $[n]$,
is
$\frac{n!}{(m!)^l(n-lm)!}$.
Given the $\{A_i\}_{i=1}^l$, the number of ways to choose the sets $\{D_j\}_{j=1}^k$ so that $\{D_j\}_{j=1}^k=\{A_i\}_{i=1}^l$ is
equal to the Stirling number $\stirling kl$, the number of ways to partition a set of size $k$ into
$l$ nonempty parts.
From these facts along with  \eqref{Cnmkthmoment} and \eqref{ldistinct},
we conclude that for $n\ge mk$,
\begin{equation}\label{finaltouch}
\begin{aligned}
&E_n(C^{(n)}_m)^k=\sum_{l=1}^k\big(\frac{(n-lm)!((m-1)!)^l}{n!}\big)\big(\frac{n!}{(m!)^l(n-lm)!}\big)\stirling kl=\\
&\sum_{l=1}^k\frac1{m^l}\stirling kl=T_k(\frac1m),
\end{aligned}
\end{equation}
proving \eqref{momentmethodCm} in the case $\theta=1$.

For the case of general $\theta$, we note that the only change that must be made in the above proof is
in \eqref{ldistinct}. Recalling that $s(n,k)$ denotes the number of permutations in $S_n$ with $k$ cycles,
we have
\begin{equation}\label{ldistincttheta}
\begin{aligned}
&E_{n;\theta}\prod_{j=1}^k1_{D_j}=\frac{((m-1)!)^l\sum_{k=1}^{n-ml}s(n-lm,k)\theta^{k+l}}{\theta^{(n)}}=
\frac{\theta^{(n-ml)}((m-1)!)^l}{\theta^{(n)}}\theta^l\sim\\
&\frac{(n-ml)!((m-1)!)^l}{n!}\theta^l=\theta^lE_n\prod_{j=1}^k1_{D_j},\ \text{as}\ n\to\infty.
\end{aligned}
\end{equation}
Thus, instead of \eqref{finaltouch}, we have
$$
E_n(C^{(n)}_m)^k\sim\sum_{l=1}^k(\frac\theta m)^l\stirling kl=T_k(\frac\theta m), \ \text{as}\ n\to\infty.
$$

We now turn to \eqref{momentmethodCms}. The method of proof is simply the natural
extension of the one used to prove \eqref{momentmethodCm}; thus, since
the notation is cumbersome we will suffice with illustrating the method by proving that
\begin{equation}\label{forfinalproof}
\lim_{n\to\infty}E_n(C^{(n)}_{m_1})^{k_1}(C^{(n)}_{m_2})^{k_2}=T_{k_1}(\frac1{m_1})T_{k_2}(\frac1{m_2}).
\end{equation}
Let $n\ge m_1k_1+m_2k_2$.
By \eqref{Cnrep}, we have
\begin{equation}\label{productexp}
E_nC^{(n)}_{m_1})^{k_1}(C^{(n)}_{m_2})^{k_2}=\sum_{\stackrel{D^{1}_j\subset[n],|D^{1}_j|=m_1}{j\in[k_1]}}
\sum_{\stackrel{D^{2}_j\subset[n],|D^{2}_j|=m_2}{j\in[k_2]}}E_n\prod_{j=1}^{k_1}1_{D^{1}_j}
\prod_{j=1}^{k_2}1_{D^{2}_j}.
\end{equation}

Now $E_n\prod_{j=1}^{k_1}1_{D^{1}_j}
\prod_{j=1}^{k_2}1_{D^{2}_j}\neq0$
if and only if for some $l_1\in[k_1]$ and some $l_2\in[k_2]$, there exist
disjoint sets $\{A^{1}_i\}_{i=1}^{l_1}$, $\{A^{2}_i\}_{i=1}^{l_2}$   such that $\{D^{r}_j\}_{j=1}^{k_r}=
\{A^{r}_i\}_{i=1}^{l_r},\ r=1,2$.
If this is the case, then
\begin{equation}\label{l12distinct}
E_n\prod_{j=1}^{k_1}1_{D^{1}_j}
\prod_{j=1}^{k_2}1_{D^{2}_j}=\frac{(n-l_1m_1-l_2m_2)!((m_1-1)!)^{l_1}(m_2-1)!)^{l_2}}{n!}.
\end{equation}
The number of ways to construct  disjoint, ordered sets  $\{A^1_i\}_{i=1}^l, \{A^2_i\}_{i=1}^l$, with
 the $A^1_i$ each consisting of $m_1$ elements from $[n]$ and the $A^2_i$ each consisting of $m_2$ elements from $[n]$, is
$\frac{n!}{(m_1!)^{l_1}(m_2!)^{l_2}(n-l_1m_1-l_2m_2)!}$.
Given $\{A^1_i\}_{i=1}^l, \{A^2_i\}_{i=1}^l$,
the number of ways to choose the ordered sets $\{D^1_j\}_{j=1}^{k_1}, \{D^2_j\}_{j=1}^{k_2}$ so that $\{D^1_j\}_{j=1}^{k_1}=
\{A^1_i\}_{i=1}^{l_1}$ and  $\{D^2_j\}_{j=1}^{k_2}=
\{A^2_i\}_{i=1}^{l_2}$   is
equal to   $\stirling {k_1}{l_1}\stirling{k_2}{l_2}$.
We have
$$
\frac{(n-l_1m_1-l_2m_2)!((m_1-1)!)^{l_1}(m_2-1)!)^{l_2}}{n!}\frac{n!}{(m_1!)^{l_1}(m_2!)^{l_2}(n-l_1m_1-l_2m_2)!}=
\frac1{m_1^{l_1}m_2^{l_2}}.
$$
From these fact along with \eqref{productexp} and \eqref{l12distinct}, we conclude that for $n\ge m_1k_1+m_2k_2$,
$$
E_nC^{(n)}_{m_1})^{k_1}(C^{(n)}_{m_2})^{k_2}=
\sum_{l_2=1}^{k_2}\sum_{l_1=1}^{k_1} \frac1{m_1^{l_1}m_2^{l_2}}\stirling {k_1}{l_1}\stirling{k_2}{l_2}=T_{k_1}(\frac1{m_1})
T_{k_2}(\frac1{m_2}),
$$
 proving   \eqref{forfinalproof}.\hfill $\square$

\section{Proofs of Theorems \ref{1} and \ref{2}}\label{Proofmulti}
\noindent \it Proof of Theorem \ref{1}.\rm\
Since $\mathcal{E}^{(n)}_m$ converges weakly to $\mathcal{E}_m$, it follows from the discussion in the first paragraph of section \ref{ProofC} that it suffices to show that
\begin{equation}
\lim_{n\to\infty}E_n(\mathcal{E}^{(n)}_m)^k=v_{k;m}.
\end{equation}
Let $n\ge km$.
For $D\subset[n]$, let $1_D(\sigma)$ equal 1 or 0 according to whether or not $\sigma\in S_n$ induces an embedded
permutation on $D$.
Then we have
\begin{equation}\label{repembedded}
\mathcal{E}^{(n)}_m(\omega)=\sum_{\stackrel{D\subset[n]}{|D|=m}}1_D(\omega),
\end{equation}
and
\begin{equation}\label{Ekthmom}
E_n(\mathcal{E}^{(n)}_m)^k=\sum_{\stackrel{D_j\subset[n],|D_j|=m}{j\in[k]}}E_n\prod_{j=1}^k1_{D_j}.
\end{equation}
There is a one-to-one correspondence between collections $\{D_j\}_{j=1}^k$, satisfying $D_j\subset[n]$ and $|D_j|=m$,
and collections $\{A_I\}_{I\subset[k]}$ of disjoint sets satisfying $A_I\subset[n]$
and satisfying
\begin{equation}\label{theequation}
\sum_{I:i\in I}l_I=m,\ \text{for all}\ i\in[k],
\end{equation}
where
\begin{equation}\label{thels}
l_I=|A_I|, \ I\subset[k].
\end{equation}
The correspondence is through the formula
\begin{equation}\label{theAs}
A_I=\big(\cap_{i\in I}D_i\big)\cap\big(\cap_{i\in[k]-I}([n]-D_i)\big),\ I\subset[k].
\end{equation}

Now $\prod_{j=1}^k1_{D_j}(\omega)=1$ if and only if
for all $I\subset[k]$,
$\sigma$ induces an embedded permutation on $A_I$.
Thus,
we have
\begin{equation}\label{expprodD}
E_n\prod_{j=1}^k1_{D_j}=\frac{(n-\sum_{I\subset[k]}l_I)!\prod_{I\subset[k]}l_I!}{n!}.
\end{equation}
Given the values $l_I=|A_I|$, $I\subset[k]$, the number of ways to construct the disjoint sets $\{A_I\}_{I\subset [k]}$
is $\frac{n!}{(n-\sum_{I\subset[k]}l_I)!\prod_{I\subset[k]}l_I!}$. Using this with
\eqref{Ekthmom}-\eqref{expprodD}, it follows that
$E_n(\mathcal{E}^{(n)}_m)^k$ equals the number of solutions
$\{l_I\}_{I\subset[k]}$ to \eqref{theequation}.

To complete the proof, we will
 show that the number of solutions to \eqref{theequation} is $v_{k;m}$.
Consider the set $\cup_{j=1}^kD_j\subset[n]$. Label the elements of this set by $\{x_i\}_{i=1}^r$.
Of course, $m\le r\le km$. Now construct the sets $\{\Gamma_i\}_{i=1}^r$ by
$\Gamma_i=\{j:x_i\in D_j\}$. By construction,
the sets $\{\Gamma_i\}_{i=1}^r$ form an $m$-cover of $[k]$ (of order $r$).
There is a one-to-one correspondence between solutions to \eqref{theequation}
and $m$ covers of $[k]$; indeed, $l_I=|\{i\in [k]:\Gamma_i=I\}|$.\hfill $\square$

\bigskip

\noindent \it Proof of Theorem \ref{2}.\rm\
We use \eqref{Einfty} to calculate
$V_m(x)=Ee^{x\mathcal{E}_m}$. We have
\begin{equation}
\begin{aligned}
&V_m(x)=E\exp(x\mathcal{E}_m)=E\exp\big(x\sum_{\stackrel{(l_1,\ldots, l_m):\sum_{j=1}^mjl_j=m}{l_j\le Z_{\frac1j}, j\in[m]}}
\prod_{j=1}^m\binom{Z_\frac1j}{l_j}\big)=\\
&\sum_{u_1\ge0,\cdots, u_m\ge0}\prod_{j=1}^m(\frac1j)^{u_j}\frac1{u_j!}e^{-\frac1j}
\exp\big(x\sum_{\stackrel{(l_1,\ldots, l_m):\sum_{j=1}^mjl_j=m}{l_j\le u_j, j\in[m]}}\prod_{j=1}^m\binom{u_j}{l_j}\big).
\end{aligned}
\end{equation}
Thus, to complete the proof, we only need to show that
\begin{equation}\label{lastformula}
\sum_{\stackrel{(l_1,\ldots, l_m):\sum_{j=1}^mjl_j=m}{l_j\le u_j, j\in[m]}}\prod_{j=1}^m\binom{u_j}{l_j}=\gamma_m(u),
\end{equation}
where
\begin{equation}\label{gammadef}
\gamma_m(u)=[z^m]\prod_{j=1}^m(1+z^j)^{u_j}.
\end{equation}
Expanding with the binomial formula, we have
\begin{equation}\label{binomexp}
\prod_{j=1}^m(1+z^j)^{u_j}=\prod_{j=1}^m\big(\sum_{l_j=0}^{u_j}\binom{u_j}{l_j}z^{jl_j}\big).
\end{equation}
From \eqref{binomexp} and \eqref{gammadef}, it follows that \eqref{lastformula} holds.
\hfill $\square$

\medskip
\noindent\bf Acknowlegement.\rm\ The author thanks Ron Holzman and Roy Meshulam for a discussion and references with regard to multisets, and Ron Holzman for the reference \cite{EFG}.

\end{document}